\documentclass[conference]{IEEEtran}
\IEEEoverridecommandlockouts
\usepackage{cite}
\usepackage{amsmath,amssymb,amsfonts}
\usepackage{algorithmic}
\usepackage{graphicx}
\usepackage{textcomp}
\usepackage{xcolor}
\def\BibTeX{{\rm B\kern-.05em{\sc i\kern-.025em b}\kern-.08em
    T\kern-.1667em\lower.7ex\hbox{E}\kern-.125emX}}

\usepackage{amsmath}
\usepackage{mathtools}
\newcommand{\Span}[1]{\mathrm{span}\{#1\}}
\newcommand{\Dim}[1]{\mathrm{dim}(#1)}
\newcommand{\Cauchy}[1]{\mathcal{C}(#1)}
\newcommand{\D}{\mathrm{d}}
\newcommand{\Lie}{\mathrm{L}}

\newcommand{\XxUlu}{\mathcal{X}\times\mathcal{U}_{[0,l_u]}}

\newcounter{DefThmLemCorRem}
\newtheorem{definition}[DefThmLemCorRem]{Definition}
\newtheorem{theorem}[DefThmLemCorRem]{Theorem}

\usepackage{xcolor}

\usepackage{enumitem}
\begin{document}
%%%%%%%%%%%%%%%%%%%%%%%%%%%%%%%%%%%%%%%%%%%%%%%%%%%%%%%%%%%%%%%%%%%%%%%%%%%%
%\bstctlcite{IEEEexample:BSTcontrol}
%%%%%%%%%%%%%%%%%%%%%%%%%%%%%%%%%%%%%%%%%%%%%%%%%%%%%%%%%%%%%%%%%%%%%%%%%%%%

\title{A Triangular Normal Form for $x$-Flat Control-Affine Two-Input Systems\\
\thanks{This research was funded in whole, or in part, by the Austrian Science Fund (FWF) P36473. For the purpose of open access, the author has applied a CC BY public copyright licence to any Author Accepted Manuscript version arising from this submission.}
}

\author{\IEEEauthorblockN{1\textsuperscript{st} Conrad Gstöttner}
\IEEEauthorblockA{\textit{Institute of Automatic Control} \\
\textit{and Control Systems Technology}\\
\textit{JKU Linz, Austria}\\
https://orcid.org/0000-0003-2107-3009}
\and
\IEEEauthorblockN{2\textsuperscript{nd} Bernd Kolar}
\IEEEauthorblockA{\textit{Institute of Automatic Control} \\
\textit{and Control Systems Technology}\\
\textit{JKU Linz, Austria}\\
https://orcid.org/0000-0001-9710-8445}
\and
\IEEEauthorblockN{3\textsuperscript{rd} Markus Schöberl}
\IEEEauthorblockA{\textit{Institute of Automatic Control} \\
\textit{and Control Systems Technology}\\
\textit{JKU Linz, Austria}\\
https://orcid.org/0000-0001-5539-7015}
}

\maketitle

\begin{abstract}
This paper is devoted to normal forms for $x$-flat control-affine systems with two inputs. We propose a general triangular normal form which contains several other normal forms discussed in the literature as special cases. We derive conditions under which a system with given $x$-flat output can be transformed into the proposed triangular form. Based on the triangular form we motivate a simple algorithm for identifying candidates for flat outputs.
\end{abstract}


\section{Introduction}
	The concept of differential flatness has been introduced by Fliess, L\'evine, Martin and Rouchon in \cite{FliessLevineMartinRouchon:1992,FliessLevineMartinRouchon:1992-2}, and has attracted a lot of interest in the control systems community. The property of a system to be differentially flat (or just ``flat'' for short) allows for a systematic solution of feed-forward and feedback control problems, see e.g. \cite{FliessLevineMartinRouchon:1995,FliessLevineMartinRouchon:1999,Rudolph:2021,Sira-RamirezAgrawal:2004}. Roughly speaking, a nonlinear control system of the form
	\begin{align}\label{eq:sys}
		\begin{aligned}
			\dot{x}&=f(x,u)
		\end{aligned}
	\end{align}
	with $\Dim{x}=n$ states and $\Dim{u}=m$ inputs is flat if there exists an $m$-dimensional (fictitious) output 
	\begin{align}\label{eq:y}
		\begin{aligned}
			y&=\varphi(x,u,\dot{u},\ldots,u^{(q)})\,,
		\end{aligned}
	\end{align}
	$u^{(q)}$ denoting the $q$-th time derivative of the input, such that the state and the input of the system can locally be expressed as functions of this output and a finite number of its time derivatives, i.e., 
	\vspace{-1ex}
	\begin{align*}
		\begin{aligned}
			x&=F_x(y,\dot{y},\ldots,y^{(r-1)})\\
			u&=F_u(y,\dot{y},\ldots,y^{(r)})\,.
		\end{aligned}
	\end{align*}
	Such a (fictitious) output \eqref{eq:y} is called a flat output of the system \eqref{eq:sys}. The computation of flat outputs is known to be a difficult problem. There do not exist easily verifiable necessary and sufficient conditions, except for certain classes of systems, including two-input drift-less systems, see \cite{MartinRouchon:1994}, systems which are linearizable by a one-fold prolongation of a suitably chosen control, see \cite{NicolauRespondek:2017}, two-input systems which are linearizable by a two-fold prolongation of a suitably chosen control, see \cite{NicolauRespondek:2016-2,GstottnerKolarSchoberl:2021-2,GstottnerKolarSchoberl:2021-3}. Necessary and sufficient conditions for $(x,u)$-flatness of control-affine systems with two inputs and four states can be found in \cite{Pomet:1997}. In the references \cite{BououdenBoutatZhengBarbotKratz:2011,SilveiraPereiraRouchon:2015,GstottnerKolarSchoberl:2021,GstottnerKolarSchoberl:2022,SchoberlSchlacher:2014}, structurally flat triangular forms have been used for computing flat outputs. Such triangular forms have proven to be very useful in this problem. E.g. in \cite{BououdenBoutatZhengBarbotKratz:2011,SilveiraPereiraRouchon:2015,GstottnerKolarSchoberl:2021,GstottnerKolarSchoberl:2022} triangular forms have been geometrically characterized, which provides easily applicable sufficient conditions for flatness. In \cite{SchoberlSchlacher:2014} an implicit triangular form together with a constructive scheme for successively transforming a system into the proposed triangular form is presented.
	
	In this paper, we propose a very general structurally flat triangular form for a particular class of flat systems, namely control-affine two input systems
	\begin{align}\label{eq:two-input affine intro}
		\begin{aligned}
			\dot{x}&=f(x)+g_1(x)u^1+g_2(x)u^2
		\end{aligned}
	\end{align}
	which are $x$-flat, i.e., which possess a flat output $y=(\varphi^1(x),\varphi^2(x))$ that depends on the system's state only. We derive conditions under which a system \eqref{eq:two-input affine intro} with given $x$-flat output can be transformed into the proposed triangular form. Although deriving a complete characterization of this triangular form does not seem to be a tractable problem, it is still useful for the problem of computing flat outputs. Indeed, based on the proposed triangular form we motivate a simple algorithm for identifying candidates for flat outputs and the triangular form may be used as a basis for formulating more sophisticated algorithms as well.
	
	The paper is structured as follows. In Section \ref{sec:notation and flatness} we introduce some notation and provide a definition for flatness. The normal form and its relation to some normal forms discussed in the literature are presented in Section \ref{sec:normal form}. A simple algorithm for identifying candidates for flat outputs is presented in Section \ref{sec:identifying candidates}. In Section \ref{sec:example} the results are illustrated by means of an example and proofs are provided in Section \ref{sec:proofs}.
	
\section{Notation and Flatness}\label{sec:notation and flatness}
	We make use of the Einstein summation convention. 
	%Let $\mathcal{X}$ be an $n$-dimensional smooth manifold, equipped with local coordinates $x^i$, $i=1,\ldots,n$. A vector field is a section of the tangent bundle, i.e., a map $v:\mathcal{X}\rightarrow\mathcal{T}(\mathcal{X})$ such that $\tau_{\mathcal{X}}\circ v=id_{\mathcal{X}}$. In local coordinates, a vector field reads $v=v^i(x)\partial_{x^i}$. Likewise, a covector field or (differential) 1-form is a section of the cotangent bundle, i.e. a map $\omega:\mathcal{X}\rightarrow\mathcal{T}^\ast(\mathcal{X})$ which in local coordinates reads $\omega=\omega_i(x)\D x^i$. A distribution on $\mathcal{X}$ of rank $k$ is a map which assigns to each point $p\in\mathcal{X}$ a $k$-dimensional linear subspace $D_p\subset\mathcal{T}_p(\mathcal{X})$ of the tangent space at $p$. Then locally there exist $k$ vector fields $v_1,\ldots,v_k$ such that $v_{1,p},\ldots,v_{k,p}$ form a basis for $D_p$. We say that the distribution $D$ is (locally) spanned by the vector fields $v_1,\ldots,v_k$, which form a (local) basis for $D$, i.e., $D=\Span{v_1,\ldots,v_k}$, with the span over the ring $C^\infty(\mathcal{X})$ of smooth functions $h:\mathcal{X}\rightarrow\mathbb{R}$. Codistributions are defined analogously. 
	The $k$-fold Lie derivative of a function $\varphi$ along a vector field $v$ is denoted by $\Lie_{v}^k\varphi$ and we define the Lie derivative of a codistribution $P=\Span{\omega^1,\ldots,\omega^p}$ by $\Lie_v P=\Span{\Lie_v\omega^1,\ldots,\Lie_v\omega^p}$.\footnote{Note that $\Lie_{v}P$ on its own is meaningless as it depends on the particular choice of generators $\omega^l$, but e.g. $P+\Lie_{v}P$ is a well defined codistribution.} The Lie bracket of two vector fields $v$ and $w$ is denoted by $[v,w]$. Given a distribution $D$, its involutive closure is denoted by $\overline D$ and we write $\Cauchy{D}$ for its characteristic distribution, which is spanned by all vector fields $v\in D$ that satisfy $[v,D]\subset D$. Characteristic distributions are always involutive. By $\partial_xh$ we denote the $m\times n$ Jacobian matrix of $h=(h^1,\ldots,h^m)$ with respect to $x=(x^1,\ldots,x^n)$. The symbols $\subset$ and $\supset$ are used in the sense that they also include equality. 
	We write $h_{[\alpha]}$ for the $\alpha$-th time derivative of a function $h$. When $h$ consists of several components, i.e., $h=(h^1,\ldots,h^m)$, then $h_{[\alpha]}=(h^1_{[\alpha]},\ldots,h^m_{[\alpha]})$. To keep expressions involving time derivatives of different orders short we use multi-indices. Let $A=(a_1,\ldots,a_m)$ be a multi-index. Then $h_{[A]}=(h^1_{[a_1]},\ldots,h^m_{[a_m]})$ and $h_{[0,A]}=(h^1_{[0,a_1]},\ldots,h^m_{[0,a_m]})$ where $h^j_{[0,a_j]}=(h^j_{[0]},\ldots,h^j_{[a_j]})$. Addition and subtraction of a multi-index with an integer $c$ is defined by $A\pm c=(a_1\pm c,\ldots,a_m\pm c)$. %Let $A=(a^1,\ldots,a^m)$ and $B=(b^1,\ldots,b^m)$ be two multi-indices with $a^j\leq b^j$, $j=1,\ldots,m$, which we abbreviate by $A\leq B$. Then $h_{[A]}=(h^1_{[a^1]},\ldots,h^m_{[a^m]})$, $h_{[0,A]}=(h^1_{[0,a^1]},\ldots,h^m_{[0,a^m]})$, $h_{[A,B]}=(h^1_{[a^1,b^1]},\ldots,h^m_{[a^m,b^m]})$, where $h^j_{[a^j,b^j]}=(h^j_{[a^j]},\ldots,h^j_{[b^j]})$. We define $h^j_{[\alpha^j,\beta^j]}$ to be empty when $\alpha^j>\beta^j$. Addition and subtraction of multi-indices is done componentwise and we define the addition and subtraction of a multi-index $A$ with an integer $c$ by $A\pm c=(a^1\pm c,\ldots,a^m\pm c)$. 
	Furthermore, we define $\#A=\sum_{j=1}^{m}a_j$. Throughout, all functions and (co)vector fields are assumed to be smooth and all (co)distributions are assumed to have locally constant rank, we consider generic points only. We call two systems $\dot{x}=f(x,u)$ and $\dot{\bar{x}}=\bar{f}(\bar{x},\bar{u})$ static feedback equivalent if they are locally equivalent via a state transformation $\bar{x}=\Phi_x(x)$ and a static feedback transformation $\bar{u}=\Phi_u(x,u)$, i.e., $(f^l(x,u)\partial_{x^l}\Phi^i_x(x))\circ\hat{\Phi}(\bar{x},\bar{u})=\bar{f}^i(\bar{x},\bar{u})$ where $\hat{\Phi}$ denotes the inverse of the transformation $(\bar{x},\bar{u})=(\Phi_x(x),\Phi_u(x,u))$.

	Like e.g. in \cite{KolarSchoberlSchlacher:2016-3}, we use a finite-dimensional differential-geometric framework with an extended state and input manifold $\XxUlu$ with coordinates $(x,u,u_{[1]},\ldots,u_{[l_u]})$, where $u_{[\alpha]}$ denotes the $\alpha$-th time derivative of the input $u$ and $l_u$ is some sufficiently large but finite integer. The time derivative $h_{[1]}$ of a function $h(x,u,u_{[1]},\ldots,u_{[l_u-1]})$ along trajectories of \eqref{eq:sys} is then given by the Lie derivative $\Lie_{f_u}h$ along the vector field
	\vspace{-1ex}
	\begin{align}\label{eq:fu}
		\begin{aligned}
			f_u&=f^i(x,u)\partial_{x^i}+\sum_{\alpha=0}^{l_u-1}u^j_{[\alpha+1]}\partial_{u^j_{[\alpha]}}\,.
		\end{aligned}
	\end{align}
	Throughout, we assume that $l_u$ is chosen large enough such that $f_u$ acts as time derivative on all functions considered. Within this framework, flatness can be defined as follows.
	\begin{definition}\label{def:flatness}
		The system \eqref{eq:sys} is called flat if there exists an $m$-tuple of smooth functions
		\begin{align}\label{eq:flat output}
			\varphi^j(x,u,u_{[1]},\ldots,u_{[q]}),\quad j=1,\ldots,m
		\end{align}
		defined on $\XxUlu$ and smooth functions $F_x^i$, $F_u^j$ such that locally the conditions
		\begin{align}\label{eq:param}
			\begin{aligned}
				x^i&=F_x^i(\varphi_{[0,R-1]}),\quad i=1,\ldots,n\\
				u^j&=F_u^j(\varphi_{[0,R]}),\quad j=1,\ldots,m
			\end{aligned}
		\end{align}
		with some multi-index $R=(r_1,\ldots,r_m)$ hold. The $m$-tuple \eqref{eq:flat output} is called a flat output.		
	\end{definition}
	Definition \ref{def:flatness} implies that the differentials $\D\varphi,\D\varphi_{[1]},\ldots,\D\varphi_{[\beta]}$ up to an arbitrary order $\beta$ are linearly independent, and based on that, it can be shown that locally there exist a unique minimal multi-index $R$ and unique maps $F_x$ and $F_u$ such that \eqref{eq:param} holds, see e.g. \cite{KolarSchoberlSchlacher:2016-3} for further details. The map $(F_x,F_u):\mathbb{R}^{\#R+m}\rightarrow\mathbb{R}^{n+m}$ is a submersion (it degenerates to a diffeomorphism for linearizing outputs in the sense of static feedback linearization). We denote the difference of the dimensions of the domain and the codomain of this map by $d$, i.e., $d=\#R+m-(n+m)=\#R-n$, and refer to it as the differential difference of the
	corresponding flat output. In \cite{NicolauRespondek:2016,NicolauRespondek:2017}, the number $\#R+m$ is called the differential weight of the flat output. (The differential weight of a flat output with differential difference $d$ is thus given by $n+m+d$.) A flat output is called minimal if its differential difference is the smallest among all flat outputs of the system.
\section{Triangular Normal Form}\label{sec:normal form}
	In the following, we consider nonlinear two-input control-affine systems
	\begin{align}\label{eq:two-input affine}
		\begin{aligned}
			\dot{x}&=f(x)+g_1(x)u^1+g_2(x)u^2
		\end{aligned}
	\end{align}
	with $\Dim{x}=n$. Let 
	\vspace{-1ex}
	\begin{align}\label{eq:x-flat output}
		y&=(\varphi^1(x),\varphi^2(x))
	\end{align}
	be an $x$-flat output of \eqref{eq:two-input affine} and let $K=(k_1,k_2)$ be the relative degrees of the flat output components, which are defined via
	\begin{align*}
		\begin{aligned}
			\Lie_{f_u}^{k_j-1}\varphi^j&=\varphi^j_{[k_j-1]}(x)\,,&&&\Lie_{f_u}^{k_j}\varphi^j&=\varphi^j_{[k_j]}(x,u)\,.
		\end{aligned}
	\end{align*}
	Since we are dealing with a control-affine system, the vector field $f_u$ is of the form
	\vspace{-1ex}
	\begin{align*}
		\begin{aligned}
			f_u&=\left(f^i(x)+u^1g_1^i(x)+u^2g_2^i(x)\right)\partial_{x^i}+\sum_{\alpha=0}^{l_u-1}u^j_{[\alpha+1]}\partial_{u^j_{[\alpha]}}
		\end{aligned}
	\end{align*}
	and $\varphi^j_{[k_j]}(x,u)$ are actually affine functions of $u$. Now consider the static feedback transformation $\bar{u}^1=\varphi^1_{[k_1]}(x,u)$, $\bar{u}^2=u^2$ (permute $u^1$ and $u^2$ if necessary).\footnote{The static feedback $\bar{u}^1=\varphi^1_{[k_1]}(x,u)$, $\bar{u}^2=u^2$ implies the transformation
		\vspace{-2ex}
	\begin{align*}
		\begin{aligned}
			\bar{u}^1_{[1]}&=\varphi^1_{[k_1+1]}(x,u,u_{[1]})&&&\bar{u}^2_{[1]}&=u^2_{[1]}\\[-1ex]
			&\vdotswithin{=}&&&&\vdotswithin{=}\\[-1ex]
			\bar{u}^1_{[l_u]}&=\varphi^1_{[k_1+l_u]}(x,u,u_{[1]},\ldots,u_{[l_u]})&&&\bar{u}^2_{[l_u]}&=u^2_{[l_u]}
		\end{aligned}
	\end{align*}
	for the derivatives of the inputs. In the new coordinates the vector field $f_u$ reads
	\vspace{-2ex}
	\begin{align*}
		\begin{aligned}
			f_u&=\underbrace{\bar{f}^i(x,\bar{u})\partial_{x^i}+\sum_{\alpha=0}^{l_u-1}\bar{u}^j_{[\alpha+1]}\partial_{\bar{u}^j_{[\alpha]}}}_{\bar{f}_u}+\ldots\partial_{\bar{u}^1_{[l_u]}}+\ldots\partial_{\bar{u}^2_{[l_u]}}\,,\\[-1ex]
		\end{aligned}
	\end{align*}
	where $\bar{f}=f\circ\hat{\Phi}_u$ and $\hat{\Phi}_u$ denotes the inverse of the transformation $\bar{u}^1=\varphi^1_{[k_1]}(x,u)$, $\bar{u}^2=u^2$ (which depends affinely on $\bar{u}$ when $\varphi^1_{[k_1]}(x,u)$ depends affinely on $u$). The in general non-zero components in the $\partial{\bar{u}^j_{[l_u]}}$-directions can be ignored as long as $l_u$ is chosen large enough, i.e., after applying a transformation, we can solely work with $\bar{f}_u$, which is of the same structure as $f_u$ in the original coordinates. (An additional state transformation $\bar{x}=\Phi_x(x)$ would only affect the $\partial_{x}$-components and result in the components $\bar{f}^i(\bar{x},\bar{u})=(f^l(x,u)\partial_{x^l}\Phi^i_x(x))\circ\hat{\Phi}(\bar{x},\bar{u})$ with respect to $\partial_{\bar{x}^i}$.)} Such transformations have also been used in \cite{GstottnerKolarSchoberl:2023} and \cite{GstottnerKolarSchoberl:2020} where it has been shown that by applying such a transformation, the flat output and its derivatives up to the orders $R=(r_1,r_2)$ take the form\footnote{The control-affine two-input $x$-flat systems studied here belong to the broader class of general nonlinear two-input systems $\dot{x}=f(x,u^1,u^2)$ with $(x,u)$-flat outputs, i.e., flat outputs $\varphi(x,u)$ which may depend on the state $x$ and also the input $u$, which are studied in \cite{GstottnerKolarSchoberl:2023} and \cite{GstottnerKolarSchoberl:2020}. The results therein are thus directly applicable here.}
%		\begin{align}\label{eq:phiR}
%			\begin{aligned}
%				\varphi^1&=\varphi^1(x)\\
%				&\vdotswithin{=}\\
%				\varphi^1_{[k^1-1]}&=\varphi^1_{[k^1-1]}(x)\\
%				\varphi^1_{[k^1]}&=\bar u^1\\
%				\varphi^1_{[k^1+1]}&=\bar{u}^1_{[1]}\\
%				&\vdotswithin{=}\\
%				\varphi^1_{[r_1]}&=\bar u^1_{[r^1-k^1]}
%			\end{aligned}\quad
%			\begin{aligned}
%				\varphi^2&=\varphi^2(x)\\
%				\vdotswithin{=}\\
%				\varphi^2_{[k^2-1]}&=\varphi^2_{[k^2-1]}(x)\\
%				\varphi^2_{[k^2]}&=\varphi^2_{[k^2]}(x,\bar u^1)\\
%				\varphi^2_{[k^2+1]}&=\varphi^2_{[k^2+1]}(x,\bar u^1,\bar{u}^1_{[1]})\\
%				&\vdotswithin{=}\\
%				\varphi^2_{[r_2]}&=\varphi^2_{[r^2]}(x,\bar u^1,\bar{u}^1_{[1]},\ldots,\\
%				&\hspace{5em}\bar{u}^1_{[r^2-k^2]},\bar{u}^2)\,.
%			\end{aligned}
%		\end{align}
	\begin{align}\label{eq:phiR}
		\begin{aligned}
			\varphi^1&=\varphi^1(x)&&&\varphi^2&=\varphi^2(x)\\[-1ex]
			&\vdotswithin{=}&&&&\vdotswithin{=}\\[-1ex]
			\varphi^1_{[k_1-1]}&=\varphi^1_{[k_1-1]}(x)&&&\varphi^2_{[k_2-1]}&=\varphi^2_{[k_2-1]}(x)\\
			\varphi^1_{[k_1]}&=\bar u^1&&&\varphi^2_{[k_2]}&=\varphi^2_{[k_2]}(x,\bar u^1)\\
			\varphi^1_{[k_1+1]}&=\bar{u}^1_{[1]}&&&\varphi^2_{[k_2+1]}&=\varphi^2_{[k_2+1]}(x,\bar u^1,\bar{u}^1_{[1]})\\[-1ex]
			&\vdotswithin{=}&&&&\vdotswithin{=}\\[-1ex]
			\varphi^1_{[r_1]}&=\bar u^1_{[r_1-k_1]}&&&\varphi^2_{[r_2]}&=\varphi^2_{[r_2]}(x,\bar u^1,\bar{u}^1_{[1]},\ldots,\\
			&&&&&\hspace{4em}\bar{u}^1_{[r_2-k_2]},\bar{u}^2)\,.\\[-1ex]
		\end{aligned}
	\end{align}
	\vspace{-2ex}
	
	\noindent
	It is immediate that the derivatives of $\varphi^1$ from order $k_1$ on are just derivatives of the new input $\bar u^1$. The form of the derivatives of $\varphi^2$ can be derived by exploiting the fact that by assumption $R=(r_1,r_2)$ are the minimal integers such that \eqref{eq:param} holds, from which it additionally follows that $r_1-k_1=r_2-k_2$, see e.g. \cite{GstottnerKolarSchoberl:2023} for details. Furthermore, it follows that $k_1+r_2=n$ since it must be
	possible to construct exactly $n$ independent functions of $x$ only
	from the functions \eqref{eq:phiR}. In conclusion, the relative degrees $K=(k_1,k_2)$, the state dimension $n$, the orders $R=(r_1,r_2)$ and the differential difference $d$ are related as follows (see also \cite{Gstottner:2023}, Proposition 3.19)
	\begin{align*}%\label{eq:integers properties}
		\begin{aligned}
			k_1+r_2&=n\,,\qquad k_2+r_1=n\,,\qquad%&d=r^1+r^2-n\\
			n-k_1-k_2=d\,.
		\end{aligned}
	\end{align*}
	From these properties, it also immediately follows that $R=K+d$.
	
	With any $x$-flat output \eqref{eq:x-flat output}, we associate the following sequence of codistributions\vspace{-1ex}
	\begin{align}\label{eq:P}
		\begin{aligned}
			P_{K-1}&=\Span{\D\varphi_{[0,K-1]}}\\
			P_{K}&=\Span{\D\varphi_{[0,K]}}\\[-1ex]
			&\vdotswithin{=}\\[-1ex]
			P_{K+d-1}&=P_{R-1}=\Span{\D\varphi_{[0,R-1]}}\\
			P_{K+d}&=P_{R}=\Span{\D\varphi_{[0,R]}}\,,
		\end{aligned}
	\end{align}
	as well as the following intersections with $\Span{\D x}$
	\begin{align}\label{eq:Q}
		\begin{aligned}
			Q_{K-1}&=P_{K-1}\cap\Span{\D x}\\
			Q_{K}&=P_{K}\cap\Span{\D x}\\[-1ex]
			&\vdotswithin{=}\\[-1ex]
			Q_{R-1}&=P_{R-1}\cap\Span{\D x}\,.
		\end{aligned}
	\end{align}
	Due to the linear independence of $\D\varphi,\D\varphi_{[1]},\ldots,\D\varphi_{[\beta]}$ up to an arbitrary order $\beta$ (see below Definition 1), we obviously have\vspace{-1ex}
	\begin{align*}%\label{eq:P seq}
		\begin{aligned}
			P_{K-1}\underset{2}{\subset}P_{K}\underset{2}{\subset}\ldots\underset{2}{\subset}P_{R-1}\underset{2}{\subset}P_R\,.
		\end{aligned}
	\end{align*}
	By the definition of the relative degree, we furthermore have $Q_{K-1}=P_{K-1}$ and from \eqref{eq:phiR} it is also immediate that
	\begin{align*}%\label{eq:Q seq}
		\begin{aligned}
			Q_{K-1}\underset{1}{\subset}Q_{K}\underset{1}{\subset}\ldots\underset{1}{\subset}Q_{R-1}=\Span{\D x}\,.
		\end{aligned}
	\end{align*}
	Except for $Q_{K-1}$ and $Q_{R-1}$ the remaining codistributions in \eqref{eq:Q} are not necessarily completely integrable. In case that all of them are indeed completely integrable, the system can be transformed into a triangular form compatible with the corresponding (not necessarily minimal) $x$-flat output $\varphi(x)$, as the following theorem asserts.	
	%TODO: Relation to $x$-maximal flatness as defined in \cite{NicolauLiRespondek:2014}?
	\begin{theorem}\label{thm:trig}
		The system \eqref{eq:two-input affine} with $x$-flat output $\varphi=(\varphi^1(x),\varphi^2(x))$ is static feedback equivalent to the triangular form
		\vspace{-2ex}
		\begin{align}\label{eq:trigform}
			\begin{aligned}
				\dot{z}^1&=z^2\\[-1ex]
				&\vdotswithin{=}\\[-1ex]
				\dot{z}^{k_1-1}&=z^{k_1}\\
				\dot{z}^{k_1}&=v^1
			\end{aligned}\quad
			\begin{aligned}
				\dot{z}^{k_1+1}&=z^{k_1+2}\\[-1ex]
				&\vdotswithin{=}\\[-1ex]
				\dot{z}^{k_1+k_2-1}&=z^{k_1+k_2}\\
				\dot{z}^{k_1+k_2}&=a^{k_1+k_2}(z^1,\ldots,z^{k_1+k_2+1})+\\
				&\hspace{6ex}b^{k_1+k_2}(z^1,\ldots,z^{k_1+k_2+1})v^1\\[1.5ex]
				\dot{z}^{k_1+k_2+1}&=a^{k_1+k_2+1}(z^1,\ldots,z^{k_1+k_2+2})+\\
				&\hspace{6ex}b^{k_1+k_2+1}(z^1,\ldots,z^{k_1+k_2+2})v^1\\[-1ex]
				&\vdotswithin{=}\\[-1ex]
				\dot{z}^{n-1}&=a^{n-1}(z)+b^{n-1}(z)v^1\\
				\dot{z}^{n}&=v^2\\[-3ex]
			\end{aligned}
		\end{align}
		with $\varphi=(z^1,z^{k_1+1})$ and $b^{k_1+k_2}\neq 0$, as well as $\partial_{z^{l+1}}a^l\neq 0$ or $\partial_{z^{l+1}}a^l\neq 0$ for $l=k_1+k_2,\ldots,n-1$, if and only if all the corresponding codistributions \eqref{eq:Q} are completely integrable.
	\end{theorem} 
	A proof of this theorem is provided in Section \ref{sec:proofs}. Due to the property $\partial_{z^{l+1}}a^l\neq 0$ or $\partial_{z^{l+1}}b^l\neq 0$ for $l=k_1+k_2,\ldots,n-1$, one of the functions $a^l$ or $b^l$ in each equation can always be normalized to $z^{l+1}$. Such normalizations preserve the triangular structure. The normal form \eqref{eq:trigform} contains several normal forms discussed in the literature as special cases:
	\begin{itemize}
		\item For $k_1+k_2=n$, i.e., if the sum of the relative degrees of the flat ouput components coincides with the state dimension, the form \eqref{eq:trigform} degenerates to the Brunovsk\'y normal form and $\varphi=(z^1,z^{k_1+1})$ is a linearizing output in the sense of static feedback linearization, see \cite{Brunovsky:1970,JakubczykRespondek:1980,HuntSu:1981,NijmeijervanderSchaft:1990,Isidori:1995}.
		\item When $k_1=k_2=1$ and all $a^l=0$, \eqref{eq:trigform} reduces (after normalizing all $b^l$) to the chained form %(also called Goursat normal form)
		\vspace{-1.5ex}
		\begin{align*}
			\begin{aligned}
				\dot{z}^{1}&=v^1
			\end{aligned}\quad
			\begin{aligned}
				\dot{z}^{2}&=z^{3}v^1\\[-1ex]
				&\vdotswithin{=}\\[-1ex]
				\dot{z}^{n-1}&=z^n v^1\\
				\dot{z}^{n}&=v^2\,,\\[-2ex]
			\end{aligned}
		\end{align*}		
		see e.g. \cite{MartinRouchon:1994,TilburySastry:1994,TilburyMurrayShankarSastry:1995,MurraySastry:1991,MurraySastry:1993,LiRespondek:2012}.
		\item The case $k_1=k_2=1$ with $a^l$ not necessarily zero and $\partial_{z^{l+1}}b^l\neq 0$ leads (after normalizing all $b^l$) to so-called extended chained systems, see \cite{LiXuSuChu:2013,NicolauLiRespondek:2014,SilveiraPereiraRouchon:2015,LiNicolauRespondek:2016}.
%		\item The two triangular forms characteriezed in \cite{GstottnerKolarSchoberl:2021} and \cite{GstottnerKolarSchoberl:2022}, which are obtained from the extended chained form by augmenting it with appropriate subsystems in Brunovsk\'y normal form are special cases of \eqref{eq:trigform}, however, transforming \eqref{eq:trigform} into one of these two forms involves in general more than just normalizing functions $a^l$ and $b^l$.
		\item In \cite{NicolauGstottnerRespondek:2022}, normal forms for $x$-flat two-input control-affine systems in dimension five are presented. All except for one of these normal forms are special cases of \eqref{eq:trigform} (the one exception is discussed below).
	\end{itemize}

	In general, not all of the codistributions \eqref{eq:Q} are completely integrable. An example is the following system in dimension five with $x$-flat output $(x^1,x^2)$
	\begin{align*}%\label{eq:nf13}
		\begin{aligned}
			\dot{x}^{1}&=u^1
		\end{aligned}\quad
		\begin{aligned}
			\dot{x}^{2}&=x^3+x^4u^1\\
			\dot{x}^{3}&=a(x^1,\ldots,x^4)+(-x^5+b(x^1,\ldots,x^4))u^1\\
			\dot{x}^{4}&=x^5+c(x^1,\ldots,x^4)u^1\\
			\dot{x}^{5}&=u^2\\[-1ex]
		\end{aligned}
	\end{align*}
	taken from \cite{NicolauGstottnerRespondek:2022}, see also \cite{NicolauLiRespondek:2014}. For the flat output $(x^1,x^2)$ the relative degrees are $K=(1,1)$ and the codistributions \eqref{eq:Q} follow as
	\vspace{-1ex}
	\begin{align*}
		Q_{(0,0)}&=\Span{\D x^1,\D x^2}\\
		Q_{(1,1)}&=\Span{\D x^1,\D x^2,\D x^3+u^1\D x^4}\\
		Q_{(2,2)}&=\Span{\D x^1,\D x^2,\D x^3,\D x^4}\\
		Q_{(3,3)}&=\Span{\D x^1,\D x^2,\D x^3,\D x^4,\D x^5}\,,
	\end{align*}
	with $Q_{(1,1)}$ being not completely integrable.

%	\begin{theorem}\label{thm:blocktrig}
%		Every $x$-flat system \eqref{eq:sys} is static feedback equivalent to the block triangular form \red{TODO: notation}
%		\begin{align}\label{eq:blocktrig}
%			\begin{aligned}
%				\dot{z}^1&=z^2\\
%				&\vdotswithin{=}\\
%				\dot{z}^{k^1-1}&=z^{k^1}\\
%				\dot{z}^{k^1}&=v^1
%			\end{aligned}\qquad
%			\begin{aligned}
%				\dot{z}^{k^1+1}&=z^{k^1+2}\\
%				&\vdotswithin{=}\\
%				\dot{z}^{k^1+k^2-1}&=z^{k^1+k^2}\\
%				\red{\dot{z}^{k^1+k^2}}&\red{=a^{k^1+k^2}(z^1,\ldots,z^{k^1+k^2+s_1})+b^{k^1+k^2}(z^1,\ldots,z^{k^1+k^2+s_1})v^1}\\
%				\blue{\dot{z}^{k^1+k^2+1}}&\blue{=a^{k^1+k^2}(z^1,\ldots,z^{k^1+k^2+s_1+s_2})+b^{k^1+k^2+1}(z^1,\ldots,z^{k^1+k^2+s_1+s_2})v^1}\\
%				&\blue{\vdotswithin{=}}\\
%				\blue{\dot{z}^{k^1+k^2+s_1}}&\blue{=a^{k^1+k^2}(z^1,\ldots,z^{k^1+k^2+s_1+s_2})+b^{k^1+k^2}(z^1,\ldots,z^{k^1+k^2+s_1+s_2})v^1}\\
%				\green{\dot{z}^{k^1+k^2+s_1+1}}&\green{=a^{k^1+k^2}(z^1,\ldots,z^{k^1+k^2+s_1+s_2+s_3})+b^{k^1+k^2+1}(z^1,\ldots,z^{k^1+k^2+s_1+s_2+s_3})v^1}\\
%				&\green{\vdotswithin{=}}\\
%				\green{\dot{z}^{k^1+k^2+s_1+s_2}}&\green{=a^{k^1+k^2}(z^1,\ldots,z^{k^1+k^2+s_1+s_2+s_3})+b^{k^1+k^2}(z^1,\ldots,z^{k^1+k^2+s_1+s_2+s_3})v^1}\\
%				&\vdotswithin{=}\\
%				\dot{z}^{n-1}&=a^{n-1}(z)+b^{n-1}(z)v^1\\
%				\dot{z}^{n}&=v^2
%			\end{aligned}
%		\end{align}
%		where the integers $s_l$ correspond to the non-zero coranks of the inclusions in \eqref{eq:B}.
%	
%	\end{theorem} 
%	\textit{Proof.} ... \hfill$\square$\\

\section{Identifying Candidates for Flat Outputs}\label{sec:identifying candidates}
	Deriving verifiable necessary and sufficient conditions for static feedback equivalence to \eqref{eq:trigform}, as it has been done for most of the specializations listed below Theorem \ref{thm:trig} in the cited references, does not seem to be a tractable problem. However, candidates for flat outputs, i.e., candidates for the top variables $(z^1,z^{k_1+1})$ in \eqref{eq:trigform}, may be identified by simple algorithms like the one we will propose below. In fact it suffices to identify only one potential component. When one component of a flat output of a two-input system is known, a second component can be computed using e.g. the results in \cite{RathinamSluis:1996}. (In this reference, the problem of computing the missing component of a flat output of a system with $m$ inputs when $m-1$ components are known, resp. have been guessed, is addressed.) As an alternative to \cite{RathinamSluis:1996}, for two-input systems \eqref{eq:two-input affine} for which one component of an $x$-flat output is known, a second component can also be computed as follows:
	\begin{enumerate}
		\item Let $\varphi^1(x)$ be a candidate for a component of an $x$-flat output. %Differentiate $\varphi^1(x)$ until an explicit dependence on an input occurs, i.e., compute $\varphi^1_{[k^1]}(x,u)$.
		\item Introduce $\varphi^1_{[k_1]}(x,u)$ as a new input via the static feedback $v^1=\varphi^1_{[k_1]}(x,u)$, $v^2=u^2$ (permute $u^1$ and $u^2$ if necessary, note that this is actually an affine feedback). 
		\item Let $\dot{x}=\bar{f}(x)+\bar{g}_1(x)v^1+\bar{g}_2(x)v^2$ be the system obtained by applying this feedback. Prolong $(n-1)$-fold the new input $v^1$, which leads to the prolonged system
%		\begin{align}\label{eq:sys prolonged}
%			\begin{aligned}
%				\dot{x}&=\bar{f}(x)+\bar{g}_1(x)v^1+\bar{g}_2(x)v^2\\
%				\dot{v}^1&=v^1_{[1]}\,,~\dot{v}^1_{[1]}=v^1_{[2]}\,,\ldots,~\dot{v}^1_{[n-2]}=v^1_{[n-1]}
%			\end{aligned}
%		\end{align}
		\begin{align}\label{eq:sys prolonged}
			\begin{aligned}
				\dot{x}&=\bar{f}(x)+\bar{g}_1(x)v^1+\bar{g}_2(x)v^2\\
				\dot{v}^1&=v^1_{[1]}\\[-1ex]
				%\dot{v}^1_{[1]}&=v^1_{[2]}\\
				&\vdotswithin{=}\\[-1ex]
				\dot{v}^1_{[n-2]}&=v^1_{[n-1]}
			\end{aligned}
		\end{align}
		with input $(v^1_{[n-1]},v^2)$ and state $(x,v^1,v^1_{[1]},\ldots,v^1_{[n-2]})$.
		\item Let %$f_p=\bar{f}^i(x)\partial_{x^i}+v^1_{[1]}\partial_{v^1}+\ldots+v^1_{[n-2]}\partial_{v^1_{[n-3]}}$, $g_{1,p}=\partial_{v^1_{[n-2]}}$, $g_{2,p}=\bar{g}_2^i(x)\partial_{x^i}$ be the drift and control vector fields of the prolonged system \eqref{eq:sys prolonged}.
		$f_p$, $g_{1,p}$, $g_{2,p}$ be the drift and control vector fields of the prolonged system \eqref{eq:sys prolonged}. Apply the well-known test for static feedback linearizability to the prolonged system \eqref{eq:sys prolonged}, i.e., compute the linearizability distributions $D_1=\Span{g_{1,p},g_{2,p}}$, $D_i=D_{i-1}+[f_p,D_{i-1}]$ for $i\geq 2$ and check them for involutivity and their ranks, see e.g. \cite{NijmeijervanderSchaft:1990}.
%		\item Check the prolonged system \eqref{eq:sys prolonged} for static feedback linearizability, see e.g. \cite{NijmeijervanderSchaft:1990}.
		\item If the prolonged system \eqref{eq:sys prolonged} is indeed static feedback linearizable, then the linearizing outputs of the prolonged system are flat outputs of the original system.% (linearizing outputs of \eqref{eq:sys prolonged} can be computed systematically from the distributions $D_i$, see e.g. \cite{NijmeijervanderSchaft:1990}).
		%\footnote{Linearizing outputs of \eqref{eq:sys prolonged} are functions of its state $(x,v^1,v^1_{[1]},\ldots,v^1_{[n-2]})$ and not neceassarily of the original system's state $x$ only. It is thus not guaranteed that the flat outputs obtained this way are $x$-flat outputs of the original system.}
		\label{it:5}
		\item If the prolonged system \eqref{eq:sys prolonged} is not static feedback linearizable, then it can be concluded that the candidate $\varphi^1(x)$ is not part of an $x$-flat output of the system and the procedure may be repeated with a different candidate.\label{it:6}
	\end{enumerate}
	The claims made in \ref{it:5}) and \ref{it:6}) can be proven based on the structure \eqref{eq:phiR}. See also \cite{GstottnerKolarSchoberl:2023}, where a similar result for $(x,u)$-flat two-input systems $\dot{x}=f(x,u^1,u^2)$ is derived.
	
	With these preliminaries at hand, let us now state a simple procedure for identifying flat output candidates, resp. candidates for one component of a flat output. Let $f=f^i(x)\partial_{x^i}$, $g_1=g_1^i(x)\partial_{x^i}$ and $g_2=g_2^i(x)\partial_{x^i}$ be the drift and control vector fields of system \eqref{eq:two-input affine}. Starting from the distribution $D_1=\Span{g_1,g_2}$, a sequence 
	\begin{align}\label{eq:distributions}
		\begin{aligned}
			D_1\subset D_2\subset \ldots \subset D_i\subset \ldots\subset \mathcal{T}(\mathcal{X})
		\end{aligned}
	\end{align}
	of distributions is computed as follows:
	\begin{itemize}
		\item If $D_i$ is involutive, then $D_{i+1}=D_i+[f,D_i]$.
		\item If $D_i$ is non-involutive, then $D_{i+1}=D_i+[D_i,D_i]$.
	\end{itemize}
	These distributions are feedback invariant and we can assume that $D_s=\mathcal{T}(\mathcal{X})$ for some integer $s$, otherwise the system is not accessible and thus certainly not flat. Let $s$ be the smallest integer such that $D_s=\mathcal{T}(\mathcal{X})$, then $D_{s-1}\underset{p}{\subset}\mathcal{T}(\mathcal{X})$ with $p\geq 1$ and $D_{s-1}$ can either be involutive or non-involutive.
	\begin{enumerate}[label=(\alph*)]
		\item If $D_{s-1}$ is involutive, then there exist $p$ functions $h^1(x),\ldots,h^p(x)$ such that $D_{s-1}^\perp=\Span{\D h^1,\ldots,\D h^p}$. Choose a function $\psi:\mathbb{R}^p\rightarrow\mathbb{R}$ and take $\varphi^1(x)=\psi (h^1(x),\ldots,h^p(x))$ as candidate for a flat output component.\label{it:a}
		\item If $D_{s-1}$ is non-involutive, then replace $D_{s-1}$ by its characteristic distribution $\Cauchy{D_{s-1}}$ and proceed as in (a).\label{it:b}
	\end{enumerate}
	In item \ref{it:a}, when $p=1$, only one candidate  has to be checked, but for $p\geq 2$, there are in fact infinitely many ways to combine the $p$ functions, i.e., there are infinitely many non-equivalent choices for $\psi$. %Regarding item \ref{it:b}, given a distribution $D$, its characteristic distribution $\Cauchy{D}$ is defined by $\Cauchy{D}\subset D$ and $[\Cauchy{D},D]\subset D$. Cauchy characteristics are always involutive and they can be computed systematically, see e.g. \cite{LiNicolauRespondek:2016,BryantChernGardnerGoldschmidtGriffiths:1991}.
	The proposed procedure is motivated by the triangular form \eqref{eq:trigform} based on the following considerations.
	\begin{itemize}
		\item In case that in \eqref{eq:trigform} we have $\partial_{z^{l+1}}b^l\neq 0$ for all $l=k_1+k_2,\ldots,n-1$, the procedure certainly succeeds:
		\begin{itemize}
			\item When in this case additionally $k_1=k_2=1$ holds, then the sequence \eqref{eq:distributions} ends with $D_s=\overline{D}_1=\mathcal{T}(\mathcal{X})$ and $\Cauchy{D_{s-1}}\underset{3}{\subset}\mathcal{T}(\mathcal{X})$. A flat output component is then indeed any function $\varphi^1$ such that $\D\varphi^1\neq 0$ and $\D\varphi^1\perp\Cauchy{D_{s-1}}$. This follows from the properties of (extended) chained systems, see \cite{LiRespondek:2012,LiNicolauRespondek:2016} for details.
			\item For $k_1\geq 2$ and/or $k_2\geq 2$, it follows that $D_{s-1}$ is involutive and its annihilator is spanned by $\D z^1$ and/or $\D z^{k_1+1}$.
		\end{itemize}
		\item When however $\partial_{z^{l+1}}b^l\neq 0$ does not hold for all $l=k_1+k_2,\ldots,n-1$, then besides obtaining new directions via $D_{i+1}=D_i+[D_i,D_i]$, there are also steps in which new directions are obtained by taking Lie brackets with the drift vector field, i.e., $D_{i+1}=D_i+[f,D_i]$. It follows that as long as at least one of the flat output components has a sufficiently large relative degree $k_j$, then this component can still be found via the proposed procedure.	
	\end{itemize}
	Finally, it should be noted that a system \eqref{eq:two-input affine} can of course have many different $x$-flat outputs for which the condition of Theorem \ref{thm:trig} is met, and these flat outputs of course need not have the same differential difference. It is thus not guaranteed that flat output components found via the proposed procedure belong to minimal flat outputs. The procedure only aims at finding a component that belongs to one of the system's flat outputs which is compatible with \eqref{eq:trigform}.
	
\section{Example}\label{sec:example}
	Consider the following model of an induction motor (already simplified by a static feedback), taken from \cite{Chiasson:1998}, see also \cite{NicolauRespondek:2016}
	\vspace{-1.5ex}
	\begin{align*}%\label{eq:motor}
		\begin{aligned}
			\dot{\theta}&=\omega\\
			\dot{\omega}&=\mu\psi_dI_q-\tfrac{\tau_L}{J}\\
			\dot{\psi}_d&=-\eta\psi_d+\eta MI_d
		\end{aligned}&\qquad
		\begin{aligned}
			\dot{\rho}&=n_p\omega+\tfrac{\eta MI_q}{\psi_d}\\
			\dot{I}_d&=v_d\\
			\dot{I}_q&=v_q\,.\\[-1ex]
		\end{aligned}
	\end{align*}
	The inputs are $v_d,v_q$. The drift and input vector fields are given by $f=\omega\partial_\theta+(\mu\psi_dI_q-\tfrac{\tau_L}{J})\partial_\omega+(-\eta\psi_d+\eta MI_d)\partial_{\psi_d}+(n_p\omega+\tfrac{\eta MI_q}{\psi_d})\partial_\rho$, $g_1=\partial_{I_d}$, $g_2=\partial_{I_q}$. Applying the procedure for identifying candidates for flat output components from the previous section, for the distributions \eqref{eq:distributions}, we obtain
	\begin{align*}
		\begin{aligned}
			D_1&=\Span{g_1,g_2}=\Span{\partial_{I_d},\partial_{I_q}}\\
			D_2&=D_1+[f,D_1]=\Span{\partial_{I_d},\partial_{I_q},\partial_{\psi_d},\psi_d^2\mu\partial_\omega+M\eta\partial_\rho}\\
			D_3&=D_2+[D_2,D_2]=\Span{\partial_{I_d},\partial_{I_q},\partial_{\psi_d},\partial_\omega,\partial_\rho}=\overline{D}_2\\
			D_4&=D_3+[f,D_3]=\mathcal{T}(\mathcal{X})
		\end{aligned}
	\end{align*}
	and since $D_{s-1}=D_3$ is involutive, item \ref{it:a} applies. We have $D_3^\perp=\Span{\D\theta}$ and thus $\varphi^1=\theta$ as candidate for a component of an $x$-flat output. It follows that $\varphi^1=\theta$ is indeed a flat output component, and a suitable second component is given by $\varphi^2=\rho$. It can be computed using e.g. the 6 step procedure from the previous section. Following the sufficiency part of the proof of Theorem \ref{thm:trig}, the system can be represented in the form \eqref{eq:trigform}. This is achieved by the state transformation $z^1=\varphi^1=\theta$, $z^2=\varphi^1_{[1]}=\omega$, $z^3=\varphi^1_{[2]}=\mu\psi_dI_q-\tfrac{\tau_L}{J}$, $z^4=\varphi^2=\rho$, $z^5=\varphi^2_{[1]}=n_p\omega+\tfrac{\eta MI_q}{\psi_d}$, $z^6=I_d$ and the static feedback $v^1=\varphi^1_{[3]}$, $v^2=u_d$, which leads to
	\begin{align*}
		\begin{aligned}
			\dot{z}^1&=z^2\\
			\dot{z}^{2}&=z^{3}\\
			\dot{z}^{3}&=v^1
		\end{aligned}\qquad
		\begin{aligned}
			\dot{z}^{4}&=z^{5}\\
			\dot{z}^{5}&=a^{5}(z^2,z^3,z^5,z^6)+b^{5}(z^2,z^3,z^5)v^1\\
			\dot{z}^{6}&=v^2\,.
		\end{aligned}
	\end{align*}

%\appendix
\section{Proof of Theorem \ref{thm:trig}}\label{sec:proofs}
	\textit{Necessity.} It is straightforward to verify that for a system in the form \eqref{eq:trigform} all the codistributions \eqref{eq:Q} computed for $\varphi=(z^1,z^{k_1+1})$ are indeed completely integrable.
	\textit{Sufficiency.} We have to show that a system \eqref{eq:two-input affine} with $x$-flat output $(\varphi^1(x),\varphi^2(x))$ for which the codistributions \eqref{eq:Q} are all completely integrable can be put into the form \eqref{eq:trigform} by means of a state transformation and a static feedback.\linebreak For that, let us first show that the required state transformation is given by\vspace{-2ex}
		\begin{align*}
			\begin{aligned}
				z^1&=\varphi^1(x)\\
				z^2&=\varphi^1_{[1]}(x)\\[-1ex]
				&\vdotswithin{=}\\[-1ex]
				z^{k_1}&=\varphi^1_{[k_1-1]}(x)	
			\end{aligned}\qquad
			\begin{aligned}			 	
				z^{k_1+1}&=\varphi^2(x)\\
				z^{k_1+2}&=\varphi^2_{[1]}(x)\\[-1ex]
				&\vdotswithin{=}\\[-1ex]
				z^{k_1+k_2}&=\varphi^2_{[k_2-1]}(x)\\
				z^{k_1+k_2+1}&=g^{k_1+k_2+1}(x)\\[-1ex]
				&\vdotswithin{=}\\[-1ex]
				z^{n}&=g^n(x)\,,\\[-1ex]
			\end{aligned}
		\end{align*}
		where the functions $g$ can be chosen as any function such that $Q_{K-1+l}=Q_{K-1}+\Span{\D g^{k_1+k_2+1},\ldots,\D g^{k_1+k_2+l}}$, $l=1,\ldots,n-k_1-k_2$. The existence of such functions $g$ is assured since by assumption all the codistributions \eqref{eq:Q} are completely integrable. In other words, we introduce new coordinates $z$ in which all the codistributions \eqref{eq:Q} are straightened out simultaneously and use $\varphi_{[0,K-1]}(x)$ as a part of these new coordinates. In these new coordinates we by construction have\vspace{-0.5ex}
		\begin{align*}
			Q_{K-1}&=\Span{\D z^1,\ldots,\D z^{k_1+k_2}}\\
			Q_{K}&=\Span{\D z^1,\ldots,\D z^{k_1+k_2+1}}\\[-1ex]
			&\vdotswithin{=}\\[-1ex]
			Q_{R-1}&=\Span{\D z^1,\ldots,\D z^{n}}\,.\\[-4ex]
		\end{align*}
		Next, introduce $v^1=\varphi^1_{[k_1]}(x,u)$ by means of a static feedback (and keep $u^2$ as the other input, which may require permuting $u^1$ and $u^2$). Note that this is actually an affine feedback. For the codistributions \eqref{eq:P} we then have
		\begin{align*}
			\begin{aligned}
				P_{K-1}&=Q_{K-1}=\Span{\D z^1,\ldots,\D z^{k_1}\!,\D z^{k_1+1},\ldots,\D z^{k_1+k_2}}\\[0.6ex]
				P_{K}&=\Span{\D z^1,\ldots,\D z^{k_1},\D v^1,\\
					&\hspace{9ex}\D z^{k_1+1},\ldots,\D z^{k_1+k_2},\D\bar\varphi^2_{[k_2]}}\\[0.6ex]
				P_{K+1}&=\Span{\D z^1,\ldots,\D z^{k_1},\D v^1,\D v^1_{[1]},\\
					&\hspace{9ex}\D z^{k_1+1},\ldots,\D z^{k_1+1},\D\bar\varphi^2_{[k_2]},\D\bar\varphi^2_{[k_2+1]}}\\[-2ex]
				&\vdotswithin{=}\\[-0.5ex]
			\end{aligned}
		\end{align*}
		where $\bar\varphi^2_{[\alpha]}$ denotes $\varphi^2_{[\alpha]}$ expressed in the new coordinates. Because of $Q_A\subset P_A$ and $\D v^1,\D v^1_{[1]},\ldots\notin Q_A$, we can construct new bases for the codistributions $P_A$, namely\vspace{-0.5ex}
		\begin{align*}
			\begin{aligned}
				P_{K-1}&=Q_{K-1}=\Span{\D z^1,\ldots,\D z^{k_1}\!,\D z^{k_1+1},\ldots,\D z^{k_1+k_2}}\\[0.6ex]
				P_{K}&=\Span{\D z^1,\ldots,\D z^{k_1},\D v^1,\\
					&\hspace{9ex}\D z^{k_1+1},\ldots,\D z^{k_1+k_2},\D z^{k_1+k_2+1}}\\[0.6ex]
				P_{K+1}&=\Span{\D z^1,\ldots,\D z^{k_1},\D v^1,\D v^1_{[1]},\\
					&\hspace{9ex}\D z^{k_1+1},\ldots,\D z^{k_1+k_2},\D z^{k_1+k_2+1},\D z^{k_1+k_2+2}}\\[-1ex]
				&\vdotswithin{=}\\
				P_{R-1}&=\Span{\D z^1,\ldots,\D z^{k_1},\D v^1,\D v^1_{[1]},\ldots,\D v^1_{[n-k_1-k_2-1]},\\
					&\hspace{9ex}\D z^{k_1+1},\ldots,\D z^{n}}\\[0.6ex]
				P_{R}&=\Span{\D z^1,\ldots,\D z^{k_1},\D v^1,\D v^1_{[1]},\ldots,\D v^1_{[n-k_1-k_2]},\\
					&\hspace{9ex}\D z^{k_1+1},\ldots,\D z^{n},\D u^2}\,,
			\end{aligned}
		\end{align*}
		where the differentials $\D\bar\varphi^2_{[k_2]},\D\bar\varphi^2_{[k_2+1]},\ldots,\D\bar\varphi^2_{[r_2-1]}$ got replaced by $\D z^{k_1+k_2+1},\D z^{k_1+k_2+2},\ldots,\D z^n$, and $\D\bar\varphi^2_{[r_2]}$ by $\D u^2$

		\noindent
		(which can be done because flatness obviously implies that $\Span{\D z,\D v^1,\D u^2}\subset P_R$). Using these bases of the codistributions $P_A$, the triangular structure of the system equations in the $(z,v)$-coordinates can be deduced as follows. Because of $P_{A+1}=P_A+\Lie_{\bar{f}_u}P_A$, with $\bar{f}_u$ denoting \eqref{eq:fu} in the newly introduced coordinates (see also footnote 2), it follows that for $i=k_1+k_2,\ldots,n-1$ the differential $\D\Lie_{\bar f_u}z^i$ can be expressed as a linear combination of the differentials $\D z^1,\ldots,\D z^{i+1},\D v^1$ and this linear combination necessarily involves $\D z^{i+1}$, but it does not involve differentials of $v^1_{[\alpha]}$, $\alpha > 0$ due to the structure of the vector field $\bar f_u$. %\footnote{Let $P=\Span{\omega^1,\ldots,\omega^k}$ be any codistribution. We define $\Lie_{f_u}P$ as the codistribution spanned by $\Lie_{f_u}\omega^1,\ldots,\Lie_{f_u}\omega^k$. Note that although $\Lie_{f_u}P$ on its own is meaningless as it depends on the particular choice of generators $\omega^l$, the codistribution $P+\Lie_{f_u}P$ is well defined.}
		This in turn is equivalent to the property that for $i=k_1+k_2,\ldots,n-1$ the time derivative $\Lie_{\bar f_u}z^i$ of $z^i$ depends on $z^1,\ldots,z^{i+1},v^1$ only and explicitly depends on $z^{i+1}$. Furthermore, due to the structure of $\bar f_u$, the derivative $\Lie_{\bar f_u}z^i$ depends affinely on $v^1$. Hence, in conclusion, there exist functions $a^i$, $b^i$ such that $\Lie_{\bar{f}_u}z^i=a^i(z^1,\ldots,z^{i+1})+b^i(z^1,\ldots,z^{i+1})v^1$ and at least $\partial_{z^{i+1}}a^i\neq 0$ or $\partial_{z^{i+1}}b^i\neq 0$. Analogously, it follows that $\Lie_{\bar{f}_u}z^n$ explicitly depends affinely on $u^2$ and can thus be normalized to $v^2$ by an affine input transformation. Together with the triangular structure of the equations for $\dot{z}^i$ for $i=1,\ldots,k_1+k_2-1$, which is an immediate consequence of the fact that we have chosen $z^i$ for $i=1,\ldots,k_1+k_2$ as successive time derivatives of flat output components, we indeed obtain the triangular form \eqref{eq:trigform}.\hfill$\square$

\bibliographystyle{IEEEtran} 
\bibliography{IEEEabrv,bibliography}

\end{document}